\documentclass[11pt]{revtex4-1}
\usepackage{amsmath, amssymb, graphicx, color}
\begin{document}
\title{Sum of digamma asymptotic error terms of an arithmetic series}
\author{Zhiqi Huang}
\affiliation{School of Physics and Astronomy, Sun Yat-sen University, 2 Daxue Road, Tangjia, Zhuhai, 519082, China}
\date{April 1st, 2023}
\def\funcS#1{\,\mathrm{S}\left(#1\right)}
\def\funcF#1#2{\,\mathrm{F}\left(#1;\, #2\right)}
\def\funcI#1#2{\,\mathrm{I}\left(#1;\, #2\right)}
\def\ccst{\frac{\ln(2\pi)-\gamma}{2}}

\begin{abstract}
  We define an $\mathrm{S}$ function as the sum of the asymptotic error terms of digamma function of an arithmetic series, $\mathrm{S}(a) \equiv \sum_{n=1}^\infty \left[\ln\frac{n}{a} - \frac{a}{2n}-\psi\left(\frac{n}{a}\right)\right]$, and show a few properties of it. Using the $\mathrm{S}$ function, we construct a real and positive $\phi$ function. Riemann hypothesis holds if $\tilde{\phi}(k)$, the complex Fourier transform of $\phi$, has only real zeros.
\end{abstract}

\maketitle

\section{The $\mathrm{S}$ function}

\begin{figure}
  \includegraphics[width=0.667\linewidth]{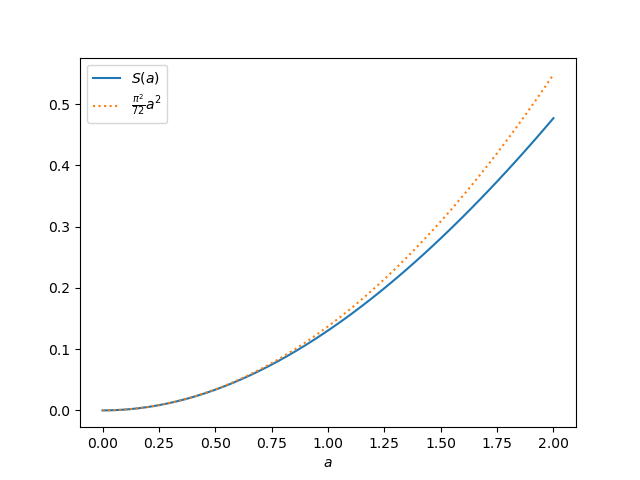}
  \caption{$\funcS{a}$ defined in Eq.~\eqref{eq:Sdef} and its asymptotic approximation $S\approx \frac{\pi^2}{72}a^2$ around $a=0$. \label{fig:Sofa}}
\end{figure}

For an argument $a>0$, we define the $\mathrm{S}$ function as 
\begin{equation}
  \funcS{a} = \sum_{n=1}^\infty \left[\ln\frac{n}{a} - \frac{a}{2n}-\psi\left(\frac{n}{a}\right)\right], \label{eq:Sdef}
\end{equation}
where $\psi(x) \equiv \frac{d\ln \Gamma(x)}{dx}$ is the digamma function. The definition~\eqref{eq:Sdef} can be extended to the complex plane $-\pi < \mathrm{arg} a<\pi$. When $a\rightarrow 0^+$, each term $\ln\frac{n}{a} - \frac{a}{2n}-\psi\left(\frac{n}{a}\right)$ vanishes as $\sim \frac{a^2}{12n^2}$. It follows then $\funcS{a}\sim \frac{\pi^2a^2}{72}$ and $\lim_{a\rightarrow 0^+}\funcS{a}= 0$. We hence define $S(0) = 0$. The function $\funcS{a}$ for $a\ge 0$, the case that will assume hereafter, is plotted in Fig.~\ref{fig:Sofa}. A very efficient numeric tool performing $\funcS{a}$ evaluation is given at \url{http://zhiqihuang.top/codes/quickS.py}. The relative accuracy of the code is  $\sim 10^{-12}$.

The $\mathrm{S}$ function satisfies a functional equation
\begin{equation}
  \funcS{a} = a \funcS{\frac{1}{a}} + \ccst (1-a) + \frac{a+1}{2}\ln a,  \label{eq:functional}
\end{equation}
where $\gamma = 0.5772156649\ldots$ is the Euler-Mascheroni constant. The proof of Eq.~\eqref{eq:functional} is rather straightforward by noticing that
\begin{equation}
  \lim_{x\rightarrow 0^+}  \left[  \int_x^\infty \frac{a dt}{\left(e^t-1\right)\left(e^{at}-1\right)} - \frac{1}{x} - \frac{1+a}{2}\ln x \right] = \funcS{a} - \ccst. \label{eq:lim}
\end{equation}
Replacing $a$ with $1/a$ on both sides of Eq.~\eqref{eq:lim} yields the desired result~\eqref{eq:functional}.

The asymptotic behavior of $\funcS{a}$ for $a\rightarrow \infty$ is $S(a)\sim \frac{a}{2}\ln a$, as Eq.~\eqref{eq:functional} implies.

Eq.~\eqref{eq:lim} can be written as an integral representation of the $\mathrm{S}$ function
\begin{equation}
  \funcS{a}  = \ccst -1 +  \frac{1+a}{2}\gamma - \int_0^\infty \left[\frac{e^{-x}}{x^2}+\frac{1-a}{2x}e^{-x}  - \frac{a}{\left(e^x-1\right)\left(e^{ax}-1\right)}\right]dx.
\end{equation}
and a more compact one
\begin{equation}
  \funcS{a} = \int_0^\infty \frac{1}{e^{\frac{x}{a}}-1}\left(\frac{1}{2}+\frac{1}{e^x-1}- \frac{1}{x}\right)dx. \label{eq:intform}
\end{equation}

Trivial manipulation of \eqref{eq:intform} gives back \eqref{eq:Sdef} (thus can be viewed as a proof of Eq.~\eqref{eq:lim}), as well as a relation between the $\mathrm{S}$ function and cosine integrals
\begin{equation}
  \funcS{a} = 2\sum_{n=1}^\infty \tau(n) \int_0^\infty\frac{\cos t}{t+\frac{2n\pi}{a}},
\end{equation}
where the divisor function $\tau(n)$ is the number of divisors of $n$, e.g., $\tau(6) = 4$, $\tau(9) = 3$.

For rational arguments, explicit integration of the right hand side of Eq.~\eqref{eq:lim} yields
\begin{equation}
  \funcS{\frac{n}{m}}  =  \ccst - \frac{1}{2m} - \frac{1}{2}\ln \frac{m}{n} + \frac{\pi}{2n} \sum_{k=1}^{n-1}\left(\frac{n}{2}- k\right)\cot\frac{mk\pi}{n} + \frac{n\pi}{2m^2} \sum_{j=1}^{m-1}\left(\frac{m}{2}- j\right)\cot\frac{nj\pi}{m}, \label{eq:rational}
\end{equation}
where the positive integers $m, n$ are relative prime.

Substituting $m=n=1$ into Eq.~\eqref{eq:rational} and applying the functional equation \eqref{eq:functional} in the infinitesimal neighborhood of $a=1$, we obtain $\funcS{1} = \frac{\ln(2\pi)-\gamma-1}{2}$ and $\mathrm{S}'(1) = \frac{1}{4}$. Other examples are $\funcS{2} = \frac{\ln(4\pi)-1-\gamma}{2}$, $\funcS{3} = \frac{\ln(6\pi)-1-\gamma}{2}+\frac{\pi}{6\sqrt{3}}$, etc.

\section{The $\phi$ function and its Fourier Transform}

\begin{figure}
  \includegraphics[width=0.667\linewidth]{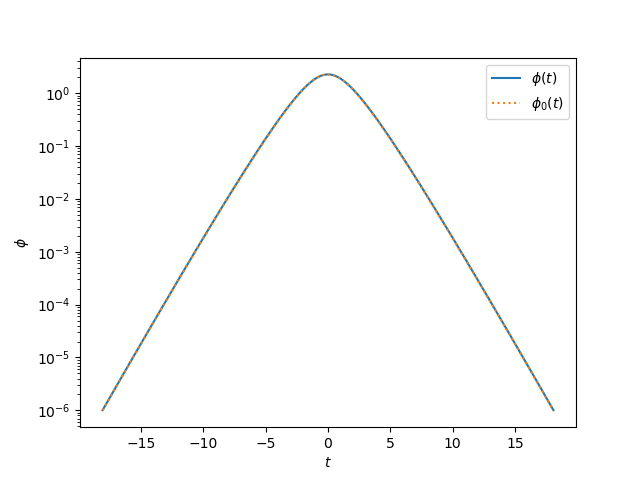}
  \caption{$\phi(t)$ defined in Eq.~\eqref{eq:phi} and its approximation $\phi_0(t)$ defined in Eq.~\eqref{eq:phi0}.\label{fig:phi}}
\end{figure}

We construct a $\phi$ function
\begin{equation}
  \phi(t) \equiv \left(t + \ccst \right) e^{-t} + e^t\funcS{e^{-2t}} + \frac{3t}{2\sinh t}, \label{eq:phi}
\end{equation}
that is positive for $t\in(-\infty, \infty)$. (At $t=0$ the $\frac{3t}{2\sinh t}$ term takes its limit value $\frac{3}{2}$.) Fig.~\ref{fig:phi} shows the $\phi$ function for real arguments.

From the functional equation~\eqref{eq:functional} we see that $\phi$ is an even function, $\phi(t) = \phi(-t)$. The Fourier transform of $\phi$,
\begin{equation}
  \tilde{\phi}(k) = \int_{-\infty}^\infty \phi(t) e^{ikt} dt, \label{eq:phiFT}
\end{equation}
is therefore also real and even for real argument $k$.

The asymptotic behavior of $\phi(t)$ as $t\rightarrow \infty$ is $\phi \sim 4te^{-t}$. The integral on the right-hand side of Eq.~\eqref{eq:phiFT} then absolutely converges for $-1<\mathrm{Im}\, k<1$. Thus, we extend the definition of $\tilde{\phi}(k)$ to the complex domain $-1<\mathrm{Im}\, k<1$ and name the right-hand side of Eq.~\eqref{eq:phiFT} as ``complex Fourier transform'' of $\phi$.

The major result of this paper is given as follows.

{{\bf Theorem}: If $z = x+iy$ ($x, y$ are real, $0<x<1$) is a nontrivial zero of Riemann zeta function $\zeta(z)$~\cite{RH}, then $\tilde{\phi}\left(2y + (1 - 2x) i\right) = 0$. }

The above theorem can be proved by taking an integral 
\begin{equation}
  \funcI{\epsilon}{z} \equiv \frac{1}{2\pi}\int_{-\pi}^{\pi} \funcF{z}{\epsilon e^{i\theta}}\funcF{1-z}{\epsilon e^{i\theta}}d\theta, \label{eq:I}
\end{equation}
where $\epsilon$ is a positive argument, and the function $\mathrm{F}$ is defined as
\begin{equation}
  \funcF{z}{w} \equiv \frac{1}{\Gamma(z)}\int_0^\infty \frac{t^{z-1}}{e^t-1+w} dt. \label{eq:F}
\end{equation}
If $0<\mathrm{Re}z<1$, the asymptotic behavior for $|w|\rightarrow 0^+$ and $-\pi < \mathrm{arg}\,w < \pi$ is
\begin{equation}
  \funcF{z}{w} = \zeta(z) + \Gamma(1-z)w^{z-1}\left(1+\frac{1+z}{2}w \right) + O(|w|),
\end{equation}
where $\zeta$ is the Riemann zeta function.

If $z = x + iy$ is a nontrivial zero of the Riemann zeta function, $\zeta(1-z)$ must also vanish, too. We then have
\begin{equation}
  \funcF{z}{w}\funcF{1-z}{w} = \frac{1}{w} \Gamma(z)\Gamma(1-z) + O(|w|^{\min(x, 1-x)}),\ \ \text{if } \ \zeta(z) = 0. 
\end{equation}
Letting $w = \epsilon e^{i\theta}$ in the integral~\eqref{eq:I}, we obtain
\begin{equation}
  \lim_{\epsilon \rightarrow 0^+} \funcI{\epsilon}{z} = 0,\ \ \text{if } \ \zeta(z) = 0. \label{eq:I0}
\end{equation}
On the other hand, we can directly compute $\funcI{\epsilon}{z}$ by substituting the definition~\eqref{eq:F} into the right-hand side of \eqref{eq:I}. The result is
\begin{equation}
  \lim_{\epsilon\rightarrow 0^+}\funcI{\epsilon}{z} = \lim_{\epsilon\rightarrow 0^+}\frac{4\sin(\pi z)}{\pi}\int_0^\infty\cos\left[\left(2y+(1-2x)i\right)t\right]h(t;\epsilon) dt, \label{eq:Ih}
\end{equation}
where
\begin{equation}
  h(t; \epsilon) \equiv  e^{-t}\left[ \int_{\epsilon e^{2t}}^\infty   \frac{ds}{\left(e^s-1\right)\left(e^{se^{-2t}}-1\right)}    -  \int_\epsilon^{\epsilon e^{2t}}  \frac{ds}{\left(e^s-1\right)\left(e^s-e^{se^{-2t}}\right)}  \right].
\end{equation}
In the limit $\epsilon\rightarrow 0^+$, $h(t; \epsilon)$ can be split into
\begin{equation}
  h(t; \epsilon) = \frac{1}{\sqrt{\epsilon }} \left[\psi \left(t+\frac{1}{2}\ln \epsilon \right) + \psi \left(-t+\frac{1}{2}\ln \epsilon \right)\right] + \phi(t), \label{eq:happ}
\end{equation}
where $\phi(t)$ is given in Eq.~\eqref{eq:phi} and
\begin{equation}
  \psi(t) \equiv e^t \int_{e^{2t}}^\infty \frac{ds}{s(e^s-1)}  - e^{-t}. \label{eq:psi}
\end{equation}
Since the right-hand side of Eq.~\eqref{eq:happ} is an even function of $t$, $\lim_{\epsilon\rightarrow 0^+}\funcI{\epsilon}{z} = 0$ then implies the complex Fourier transform of $h(t; \epsilon)$ vanishes as $\epsilon$ approaches $0^+$ and the wave number $k = 2y + (1-2x)i$. Suppose the complex Fourier transform of $\psi$ and $\phi$ are $\tilde{\psi}$ and $\tilde{\phi}$, respectively, the complex Fourier transform of $h(t; \epsilon)$ is
\begin{equation}
  \lim_{\epsilon\rightarrow 0^{+}}\tilde{h}(k; \epsilon) =  \lim_{\epsilon\rightarrow 0^{+}} \left[\tilde{\psi} \left(k\right) \epsilon^{z} + \tilde{\psi}(-k) \epsilon^{1-z}\right] + \tilde{\phi}(k)  = \tilde{\phi}(k), \ \ k = 2y+(1-2x)i. \label{eq:hphi}
\end{equation}
Combining Eqs.~\eqref{eq:I0}, \eqref{eq:Ih}  and \eqref{eq:hphi}  yields the desired conclusion.
  
\section{The path to study zeros of $\tilde{\phi}(k)$ }

The theorem in the last section implies that Riemann hypothesis (RH) holds if $\tilde{\phi}$ has only real zeros. There are many examples of positive and even functions whose complex Fourier transform has only real zeros, such as the rectangular function, the triangular function, etc. It is not easy, however, to determine whether the complex Fourier transform of a general positive and even function has only real zeros~\cite{Dimitrov11}. 

We may get some hints, however, if we consider a very good approximation  of $\phi(t)$ in a simple elementary form
\begin{equation}
  \phi_0(t) = \frac{t}{\sinh(\omega t)\cosh\left[(1-\omega)t\right]}, \label{eq:phi0}
\end{equation}
where $\omega = \frac{1}{1+\ln(2\pi)-\gamma}$. The comparison between $\phi$ and $\phi_0(t)$ is shown in Fig.~\ref{fig:phi}. The relative difference between $\phi$ and $\phi_0$ is $\sim O\left(10^{-2}\right)$ for $|t|\sim O\left(1\right)$ and drops exponentially as $|t|$ increases.

The Fourier transform of $\phi_0(t)$ is
\begin{equation}
  \tilde{\phi}_0(k) = 2\pi^2 \sum_{n=1}^\infty (-1)^{n-1} \left\{   \frac{n e^{-n\pi k/\omega}}{\omega^2 \cos\frac{n\pi(1-\omega)}{\omega}} + \frac{(n-1/2) e^{-(n-1/2)\pi k/(1-\omega)}}{(1-\omega)^2\sin\frac{(n-1/2)\pi\omega}{1-\omega}}\right\}. \label{eq:phi0k}
\end{equation}
for $\mathrm{Re}\, k > 0$. For $\mathrm{Re}\, k < 0$ we use $\tilde{\phi}_0(k) = \tilde{\phi}_0(-k)$ to evaluate $\tilde{\phi}_0(k)$.

\begin{figure}
  \includegraphics[width=0.667\linewidth]{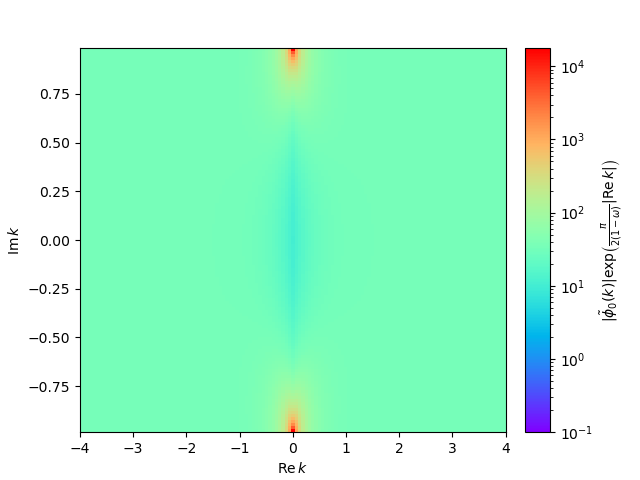}
  \caption{$\left\vert\tilde{\phi}_0(k)\right\vert \exp{\left(\frac{\pi}{2(1-\omega)} \left\vert\mathrm{Re}\,k\right\vert\right)}$ as a function of $k$ \label{fig:phi0k}}
\end{figure}

For large $\left\vert\mathrm{Re}\, k\right\vert$, the right-hand side of Eq.~\eqref{eq:phi0k} is dominated by the first $n=1$ term. When $|\mathrm{Re}\, k|$ increases,  $\left\vert\tilde{\phi}_0(k)\right\vert \exp{\left(\frac{\pi}{2(1-\omega)} \left\vert\mathrm{Re}\,k\right\vert\right)}$ exponentially converges towards a constant. Fig.~\ref{fig:phi0k} shows the numerically evaluated $\left\vert\tilde{\phi}_0(k)\right\vert \exp{\left(\frac{\pi}{2(1-\omega)} \left\vert\mathrm{Re}\,k\right\vert\right)}$ as a function of $k$.

It follows then $\tilde{\phi}_0(k)$ essentially has no zeros in the band $-1<\mathrm{Im}\,k<1$.

The accuracy of approximation can be improved, for instance, by a more complex function
\begin{equation}
  \phi_1(t) = \frac{t+\frac{\ln(2\pi)-\gamma}{4\pi}\arctan\frac{8t}{9}}{\sinh(\omega_1 t)\cosh\left[(1-\omega_1)t\right]}, \label{eq:phi1}  
\end{equation}
where $\omega_1 = \frac{1+\frac{\ln(2\pi)-\gamma}{4\pi}}{1+\ln(2\pi)-\gamma}$. The relative difference between $\phi_1(t)$ and $\phi(t)$ is $\sim 10^{-3}$ for $t\sim O(1)$ and drops exponentially as $|t|$ increases. It can also be shown that $\phi_1(t)$ also has only real zeros. It is unclear, however, whether the effort of finding better and better elementary approximations of $\phi$ can lead to a final success of proving (or disproving) that $\tilde{\phi}(k)$ has only real zeros.

\section{Conclusions and Discussion}

It has been known for more than a century that Riemann hypothesis holds if the complex Fourier transform of the positive and even function
\begin{equation}
  g(t) = \sum_{n=1}^\infty\left(2\pi n^2e^{2u}-3\right)n^2e^{5u/2 - \pi n^2e^{2u}}
\end{equation}
has only real zeros. Until so far attempts to prove that $g(t)$ has only real zeros are not very successful. (See the review article~\cite{Dimitrov11} and references therein.)

In this manuscript we have constructed another $\phi(t)$ function with the same property. While $g(t)$ contains exponentials of exponential that are in general difficult to deal with, $\phi(t)$ can be very well approximated by a simple elementary function $\phi_0(t)$ given in Eq.~\eqref{eq:phi0}. We show with straightforward calculation that $\phi_0$ does not have imaginary zeros. For better approximations with more sophisticated forms, however, it becomes more and more difficult to prove the reality of zeros of their complex Fourier transforms. The path towards a final proof of RH remains unclear.

Here we have focused on the $\mathrm{S}$ function with real arguments. Further study on the $\mathrm{S}$ function with complex variable might provide us with another path towards direct evaluation of the complex Fourier transform of $\phi$. We leave exploration along this direction as our future work.

\section{Acknowledgments}

I am grateful to my son Mr. Ningyuan Huang for his enthusiasm in RH, which motivated me, an astrophysicist, spend quite significant amount of time on studying this pure math problem. I thank my friends X. Yuan and Y. Zhou for their encouraging me to submit this manuscript to the community.


\end{document}